\documentclass[final]{siamltex}
\usepackage{a4} 
\usepackage{amsmath}
\usepackage{amssymb}
\usepackage{graphicx}
\usepackage{caption,subcaption}
\usepackage{multirow}
\usepackage{xstring}
\usepackage[utf8]{inputenc}
\usepackage{amstext,amsbsy,amsopn,eucal,enumerate}
\usepackage{tikz}
\usepackage{pgfplots}
\usepackage{tabularx}

\newcommand{\beq}{\begin{equation}}
\newcommand{\eeq}{\end{equation}}
\newcommand {\mat}      [1] {\left[\begin{array}{#1}}
\newcommand {\rix}          {\end{array}\right]}
\newcommand {\smat}      [1] {\left[\begin{smallmatrix}{#1}}
\newcommand {\srix}          {\end{smallmatrix}\right]}
\newcommand {\s}      [1] {\begin{smallmatrix}{#1}}
\newcommand {\se}          {\end{smallmatrix}}

\newtheorem{defn}{Definition}[section]

\newtheorem{lem}[defn]{Lemma}

\newtheorem{thm}[defn]{Theorem}

\title{Energy estimates and model order reduction for stochastic bilinear systems}

\author{Martin Redmann\thanks{Weierstrass Institute for Applied Analysis and Stochastics, Mohrenstrasse 39, 10117 Berlin Germany  (Email: {\tt 
martin.redmann@wias-berlin.de})}. The author~gratefully acknowledge the support from the DFG through the research unit FOR2402.}

\begin{document}

\maketitle

\begin{abstract}
In this paper, we investigate a large-scale stochastic system with bilinear drift and linear diffusion term. Such high dimensional systems appear for example when 
discretizing a stochastic partial differential equations in space. We study a particular model order reduction technique called balanced truncation (BT)
to reduce the order of spatially-discretized systems and hence reduce computational complexity. We introduce suitable Gramians to the system and prove energy estimates that can be used 
to identify states which contribute only very little to the system dynamics. When BT is applied the reduced system is obtained by removing these states from the original system. The main contribution of this paper is
an $L^2$-error bound for BT for stochastic bilinear systems. This result is new even for deterministic bilinear equations. In order to achieve it, we develop a new technique which is not 
available in the literature so far.
\end{abstract}



\pagestyle{myheadings}
\thispagestyle{plain}
\markboth{M. REDMANN}{ENERGY ESTIMATES AND MOR FOR STOCHASTIC BILINEAR SYSTEMS}


\section{Introduction}

Many phenomena in real life can be described by partial differential equations (PDEs). Famous examples are the 
motion of viscous fluids, the description of water or sound waves and the 
distribution of heat. For an accurate mathematical modeling of these real world 
applications it is often required to take random effects into account. 
Uncertainties in a PDE model can, for example, be represented by an 
additional noise term. This leads to stochastic PDEs (SPDEs) \cite{dapratozab,newspde,zabczyk,prevotroeckner}.\smallskip

It is necessary to discretize a time-dependent SPDE in space and time in order 
to solve it numerically. Discretizing in space can be 
considered as a first step. This can be done for example by spectral Galerkin \cite{galerkin, galerkinhaus, galerkinjentzen} or finite 
element methods \cite{MR1637047,MR2646102,Kruse14}. This usually leads to a high dimensional SDE. 
Solving such complex SDE systems causes large computational cost. In this context, model order reduction (MOR) is used to save computational time by
replacing large scale systems by systems of low order in which the main information of the original system should be captured. 

\subsection{Setting}\label{settingstochstabgen}

We consider a large-scale stochastic control system with bilinear drift that can be interpreted as a spatially-discretized SPDE. The corresponding noise process is an $\mathbb R^{v}$-valued L\'evy 
process $M=\left(M_1, \ldots, M_{v}\right)^T$ with mean zero and $\mathbb E \left\|M(t)\right\|^2_2=\mathbb E\left[ M^T(t)M(t)\right]<\infty$ for all $t\geq 0$. We investigate the system
     \begin{subequations}\label{controlsystemoriginal}
\begin{align}
d x(t)&=[A x(t)+ Bu(t) + \sum_{k=1}^m N_k x(t) u_k(t)]dt+ \sum_{i=1}^v H_i x(t-) dM_i(t),\label{stateeq}\\ 
y(t)&= {C} x(t),\;\;\;t\geq 0.\label{originalobserveq}
\end{align}
\end{subequations}
We assume that $A, N_k, H_i\in \mathbb R^{n\times n}$ 
($k\in\left\{1, \ldots, m\right\}$ and $i\in\left\{1, \ldots, v\right\}$), $B\in \mathbb R^{n\times m}$ and $C\in \mathbb R^{p\times n}$. Moreover, we define $x(t-):=\lim_{s\uparrow t} x(s)$.
The control $u=\left(u_1, \ldots, u_{m}\right)^T$ is assumed to be deterministic and square integrable, i.e., 
     \begin{align*}
\left\|u\right\|_{L^2_T}^2:=\int_0^T \left\|u(t)\right\|_2^2 dt<\infty
\end{align*}                            
for every $T>0$. We denote the covariance matrix of $M$ by $K=\left(k_{ij}\right)_{i, j=1, \ldots, v}$. It characterizes the covariance function of $M$, see \cite[Theorem 4.44]{zabczyk}, meaning that 
$\mathbb E[M(t)M^T(t)]=K t$. \smallskip 

We aim to replace the large scale system (\ref{controlsystemoriginal}) by a system of the same structure, but with a much smaller state dimension $r\ll n$. This reduced order model (ROM) 
is supposed be chosen, such that the corresponding output $y_r$ is close to the original one, i.e., $y_r\approx y$ in some metric.\\
In this paper, we consider balanced truncation (BT) for obtaining a ROM. It relies on defining a reachability Gramian $P$ and an observability Gramian $Q$.
These matrices are selected, such that $P$ characterizes the states in (\ref{stateeq}) and $Q$ the states in (\ref{originalobserveq}) which barely contribute to the system dynamics (see Section \ref{sec:reach} for details). 
In order to be able to remove both the unimportant states in (\ref{stateeq}) and (\ref{originalobserveq}) 
simultaneously, the first step of BT is a state space transformation 
\begin{align*}
 (A, B, C, H_i, N_k)\mapsto (\tilde A, \tilde B, \tilde C, \tilde H_i, \tilde N_k):=(SAS^{-1}, SB, CS^{-1}, SH_iS^{-1}, SN_kS^{-1}),
\end{align*}
where $S=L_QX \Sigma^{-\tfrac{1}{2}}$ and $S^{-1}=L_PY\Sigma^{-\tfrac{1}{2}}$. The ingredients of the balancing transformation are computed by the Cholesky factorizations $P=L_PL_P^T$, 
$Q=L_QL_Q^T$, and the singular value decomposition $X\Sigma Y^T=L_Q^TL_P$. This transformation does not change the output $y$ of the system, but it guarantees that the new Gramians are
diagonal and equal, i.e., $S P S^T=S^{-T}Q S^{-1}=\Sigma=\diag(\sigma_1,\ldots, \sigma_n)$ with $\sigma_1\geq \ldots \geq \sigma_n$ being the Hankel singular values (HSVs) of the system. The ROM is then obtained by 
selecting the left upper $r\times r$ blocks of $\tilde A, \tilde H_i, \tilde N_k$, the first $r$ rows of $\tilde B$ and the first $r$ columns of $\tilde C$, such that the smallest $n-r$ HSVs are removed from the system.

\paragraph{Stochastic linear case ($N_k=0$)}
BT is a method that was developed for deterministic linear systems ($H_i=N_k=0$) \cite{antoulas,moore,obiand}. Basically, two ways of extending BT to stochastic linear systems have been considered so far.
The so-called type I approach relies on defining Gramians based on generalized fundamental solutions of the system  \cite{bennerdamm, redmannbenner}. The drawback of this generalization is that only an $\mathcal H_2$-error bound 
is available \cite{redmannbenner} and an $\mathcal H_\infty$-error bound cannot be achieved \cite{bennerdammcruz, dammbennernewansatz}. To overcome this issue the type II ansatz was introduced in 
\cite{dammbennernewansatz}, where an $\mathcal H_\infty$-error bound is proved. There, a different reachability Gramian was considered which is defined as the solution to a matrix inequality. Energy estimates for linear
stochastic systems for both the type I and the type II ansatz have recently been given in \cite{redmannspa2}, such that MOR based on both approaches can be justified. As an alternative to BT, we want to refer to the 
singular perturbation approximation, where the work in \cite{fernic,spa} was extended to stochastic linear systems in \cite{redmannspa2, redSPA}.
\paragraph{Deterministic bilinear case ($H_i=0$)}
Although the bilinear term is a very weak nonlinearity, deterministic bilinear system can be seen as bridge between linear and nonlinear systems. This is because 
many deterministic nonlinear systems can be represented by bilinear systems using a so-called Carleman linearization. Applications of these equations can be found in various fields \cite{brunietal,Mohler,rugh}.
Apart from balancing related MOR techniques \cite{bennerdamm,hartmann}, various alternative methods have already been studied for the case of $H_i=0$ \cite{bennerbreiten1, bennerbreiten2, breitendamm, flaggserkan}.
We, however, only discuss the case of BT below. 
When considering BT for deterministic bilinear systems, Gramians have to be chosen properly in order to find suitable characterizations for the reachability and observability energy of the system. Concerning the choice 
of the bilinear Gramians the control components $u_k$ that are multiplied with the state $x$ in (\ref{stateeq}) are treated like white noise. Gramians according to the stochastic type I approach were, e.g., considered 
in \cite{bennerdamm, graymesko}. In both references energy estimates can be found, but they are only valid for states being in a possibly very small neighborhood of zero. Moreover, no error bound exists for this approach. 
Choosing the bilinear Gramian according to the stochastic type II approach has been considered in \cite{redmanntypeiibilinear}. To be more precise, perturbed type II Gramians were used there. This has the advantage of finding 
global energy estimates under the assumption of having bounded controls. Furthermore, an $\mathcal H_\infty$-error bound was proved in \cite{redmanntypeiibilinear}, again, assuming bounded controls. Depending on the underlying system, 
the bound on the controls can be small which is the drawback of this method.

\subsection{Outline of the paper and main result}

The work in this paper on BT for system (\ref{controlsystemoriginal}) can be interpreted as a generalization of the deterministic bilinear case, where the control $u$ is perturbed by L\'evy noise. We see 
this extension as a first step to find a bridge between stochastic linear systems and stochastic nonlinear systems to open the field of balancing related MOR to many more stochastic equations and applications.\smallskip

In this paper, the main contributions are energy estimates for the stochastic bilinear system and an error bound for BT. It is important to notice that these results are not just an extension of existing theory, they are  
new even for the deterministic bilinear case ($H_i=0$), since in contrast to \cite{redmanntypeiibilinear} no bound on the control is assumed.
To be more precise, in Section \ref{sec:reach}, we propose Gramians $P$ and $Q$ to system (\ref{controlsystemoriginal}).
We show that the reachability Gramian $P$ provides information about the degree of reachability of a state. Moreover, we establish a bound on the observation energy using the observability Gramian $Q$.
The following result on an $L^2$-error bound for BT is proved in Section \ref{errorboundsBT}. For the proof of this theorem, the existing methods in the literature cannot be applied. Hence,
we also provide a new technique to achieve this bound.
\begin{thm}\label{mainthmintro}
Let $y$ be the output of the full model (\ref{controlsystemoriginal}) with $x(0)=0$ and $y_r$ be the output of the ROM by BT with zero initial state. Then, for all $T>0$, 
it holds that \begin{align*}
  \left(\mathbb E\left\|y-y_r\right\|_{L^2_{T}}^2\right)^{\frac{1}{2}}\leq 2 (\sigma_{r+1}+\sigma_{r+2}+\ldots + \sigma_n)  \left\|u\right\|_{L^2_T}\exp\left(0.5 \left\|u\right\|_{L^2_T}^2\right),            
  \end{align*}
where $\sigma_{r+1}, \sigma_{r+2}, \ldots, \sigma_n$ are the smallest $n-r$ HSVs of system (\ref{controlsystemoriginal}).
\end{thm}\\
Theorem \ref{mainthmintro} implies that BT works well for stochastic bilinear systems if the truncated HSVs are small and the control energy is not too large. 

\section{Energy estimates}
\label{sec:reach}
BT relies on the idea to create a system (\ref{controlsystemoriginal}), in which the dominant reachable and observable states are the same. Afterwards, the unimportant states are removed to 
obtain an accurate approximation to the original model. In order to find the states that are hardly reachable and observable, a reachability Gramian $P$ and an observability Gramian $Q$ are
introduced in this section. We will see that the definitions of the Gramians are meaningful, since they lead to estimates, which allow us to find the states that barely contribute to the system dynamics. This justifies 
to balance the system based on the proposed Gramians.

\subsection{Reachability Gramian}
We introduce a reachability Gramian $P$ as a positive definite solution to 
\begin{align}\label{newgram2}
 A^T P^{-1}+P^{-1}A+\sum_{k=1}^m N^T_k P^{-1} N_k + \sum_{i, j=1}^v H_i^T P^{-1} H_j k_{i j} \leq -P^{-1}BB^T P^{-1}.
                                       \end{align}
An inequality is considered in (\ref{newgram2}), since the existence of a positive definite solution is not ensured when having an equality. The 
existence of a solution to (\ref{newgram2}) goes back to \cite{dammbennernewansatz, redmannspa2} and is given if
\begin{align}\label{stochstab}
\lambda\left(A\otimes I+I\otimes A+\sum_{k=1}^m N_k\otimes N_k+\sum_{i, j=1}^v H_i\otimes H_j k_{ij}\right)\subset \mathbb C_-, 
\end{align}
which we assume to hold througout the remainder of the paper. Here, $\lambda\left(\cdot\right)$ denotes the spectrum of a matrix. 
Condition (\ref{stochstab}) is called mean square asymptotic stability \cite{damm, staboriginal, redmannspa2}. It means that if the control components $u_k$ in the bilinear term of (\ref{stateeq}) would truely 
be white noise, then the second moment of the solution would tend to zero in the uncontrolled setting ($B=0$) if $t\rightarrow \infty$. \smallskip

Let $x(t, x_0, u)$ denote the solution to (\ref{stateeq}) for $t\geq 0$, an initial state $x_0\in\mathbb R^n$ and a control $u\in L^2_T$. We choose $(p_{k})_{k=1, \ldots, n}$ to be an 
orthonormal basis of $\mathbb R^n$ consisting of eigenvectors of $P$. We denote the corresponding eigenvalues by $(\lambda_{k})_{k=1, \ldots, n}$. For the Fourier coefficients of $x(t, 0, u)$, we obtain 
\begin{align}\nonumber
\langle x(t, 0, u), p_{k}  \rangle_2^2 &\leq \lambda_{k}\; \sum_{i=1}^n \lambda_{i}^{-1} \langle x(t, 0, u), p_{i}  
\rangle_2^2=\lambda_{k} \;\left\|\sum_{i=1}^n \lambda_{i}^{-\frac{1}{2}} \langle x(t, 0, u), p_{i}  \rangle_2 \;p_{i}\right\|_{2}^2
\\& =\lambda_{k}\; \left\|P^{-\frac{1}{2}} x(t, 0, u)\right\|_{2}^2 = \lambda_{k} \;x^T(t, 0, u) P^{-1} x(t, 0, u). \label{fourierestimate}
\end{align}
We use a shorter notation for the state below, i.e., we write $x(t)$ instead of $x(t, 0, u)$ if required. By Lemma \ref{lemstochdiff}, we have \begin{align}\label{firsteqreach}
\mathbb E\left[x^T(t) P^{-1} x(t)\right]=&2 \int_0^t \mathbb E\left[x^T(s) P^{-1} \left(A x(s)+ Bu(s) + \sum_{k=1}^m N_k x(s) u_k(s)\right)\right] ds
\\&+\int_0^t \mathbb E\left[x^T(s)\sum_{i, j=1}^v H_i^T P^{-1} H_j k_{ij} x(s)\right] ds. \nonumber
\end{align}
The bilinear term in the above equation can be bounded as follows:
\begin{align}\nonumber
&\sum_{k=1}^m 2 \int_0^t x^T(s) P^{-1} N_k x(s) u_k(s) ds=\sum_{k=1}^m 2 \int_0^t \left\langle  P^{-\frac{1}{2}} x(s)u_k(s), P^{-\frac{1}{2}} N_k 
x(s) \right\rangle_2 ds\\&\leq \sum_{k=1}^m \left( \int_0^t \left\|  P^{-\frac{1}{2}} x(s)u_k(s)\right\|_2^2ds + \int_0^t\left\|P^{-\frac{1}{2}} N_k \nonumber
x(s) \right\|^2_2 ds\right)\\&=\int_0^t  x^T(s) P^{-1} x(s) \left\|u(s)\right\|_2^2 ds + \int_0^t x^T(s) \sum_{k=1}^m N_k^TP^{-1} N_k x(s)ds.\label{maakebilinaerstoch}
\end{align}
We insert this inequality into (\ref{firsteqreach}), such that \begin{align}\label{ieneqforpinv}
&\mathbb E\left[x^T(t) P^{-1} x(t)\right]\\ \nonumber
&\leq\mathbb E\int_0^t x^T(s) (A^T P^{-1}+P^{-1} A+ \sum_{k=1}^m N_k^TP^{-1} N_k +\sum_{i, j=1}^v H_i^T P^{-1} H_j k_{ij}) x(s)ds\\ \nonumber
&\quad+\mathbb E\int_0^t 2 x^T(s) P^{-1} B u(s)ds +\int_0^t \mathbb E\left[x^T(s) P^{-1} x(s)\right] \left\|u(s)\right\|_2^2 ds.
\end{align}
We then plug in (\ref{newgram2}), which yields \begin{align*}
\mathbb E\left[x^T(t) P^{-1} x(t)\right]\leq&-\mathbb E\int_0^t x^T(s) P^{-1}BB^T P^{-1} x(s)ds\\
&+\mathbb E\int_0^t 2 x^T(s) P^{-1} B u(s)ds +\int_0^t \mathbb E\left[x^T(s) P^{-1} x(s)\right] \left\|u(s)\right\|_2^2 ds\\
=& \mathbb E \int_0^t \left\|u(s)\right\|_2^2 - \left\|B^T P^{-1} x(s)-u(s)\right\|_2^2ds\\&+\int_0^t \mathbb E\left[x^T(s) P^{-1} x(s)\right] \left\|u(s)\right\|_2^2 ds\\
\leq & \int_0^t \left\|u(s)\right\|_2^2ds+\int_0^t \mathbb E\left[x^T(s) P^{-1} x(s)\right] \left\|u(s)\right\|_2^2 ds
\end{align*}
The Gronwall inequality, see Lemma \ref{gronwall}, provides \begin{align*}
             \mathbb E\left[x^T(t) P^{-1} x(t)\right] \leq \int_0^t \left\|u(s)\right\|_2^2ds \exp\left(\int_0^t \left\|u(s)\right\|_2^2ds\right).                                               
                                                            \end{align*}
Consequently, by (\ref{fourierestimate}), we have \begin{align}\label{diffreachjaneintype2}
\sup_{t\in[0, T]}\sqrt{\mathbb E \langle x(t, 0, u), p_{k}  \rangle_2^2} \leq \lambda_{k}^{0.5} \left\|u\right\|_{L^2_T} \exp\left(0.5 \left\|u\right\|_{L^2_T}^2\right).
\end{align}
Given a state it is not possible to gain information about the corresponding energy from (\ref{diffreachjaneintype2}). However, given a bound on the energy which is not too large, let us say 
$\left\|u\right\|_{L^2_T}\leq 1$, we can conclude how much a state component contributes to the systems dynamics. If $\lambda_{k}$ is small, (\ref{diffreachjaneintype2}) implies that the Fourier coefficient
$\langle x(\cdot, 0, u), p_{k}  \rangle_2$ is close to zero on $[0, T]$ for normalized controls $u$. This means that the state variable takes only very small values in the direction of $p_{k}$ such that
hardly reachable states have a large component in the eigenspaces of $P$ belonging to the small eigenvalues. Inequality (\ref{diffreachjaneintype2}) has already been pointed out
in \cite[Remark 1]{redmanntypeiibilinear} for the case $H_i=0$. 

\subsection{Observability Gramian} 
We define the observability Gramian to be the solution to
\begin{align}\label{gengenlyapobs}
A^T Q+Q A+\sum_{k=1}^m N_k^T Q N_k +\sum_{i, j=1}^v H_i^T Q H_j k_{ij} = -C^T C.
\end{align}
Condition (\ref{stochstab}) guarantees the existence of a positive semidefinite solution \cite{redmannspa2}, but we will furthermore assume that $Q>0$ for the rest of the paper. 
Again, we use a short notation by setting $x_{x_0}(t):=x(t, x_0, u)$. In order to find a suitable estimate for $\mathbb E\left[x^T_{x_0}(t) Q x_{x_0}(t)\right]$, it is only required to replace $P^{-1}$ by $Q$ in 
(\ref{ieneqforpinv}) and take into account the additional term $x_0^T Q x_0$ that is due to the non-zero initial condition. Hence, \begin{align*}
&\mathbb E\left[x^T_{x_0}(t) Q x_{x_0}(t)\right]\\ 
&\leq\mathbb E\int_0^t x^T_{x_0}(s) (A^T Q+Q A+ \sum_{k=1}^m N_k^T Q N_k +\sum_{i, j=1}^v H_i^T Q H_j k_{ij}) x_{x_0}(s)ds\\
&\quad+\mathbb E\int_0^t 2 x^T_{x_0}(s) Q B u(s)ds +\int_0^t \mathbb E\left[x_{x_0}^T(s) Q x_{x_0}(s)\right] \left\|u(s)\right\|_2^2 ds+x_0^T Q x_0.
\end{align*} 
Applying (\ref{gengenlyapobs}) gives us \begin{align*}
\mathbb E\left[x^T_{x_0}(t) Q x_{x_0}(t)\right]\leq &-\mathbb E\int_0^t y^T(s) y(s)ds+\mathbb E\int_0^t 2 x^T_{x_0}(s) Q B u(s)ds+x_0^T Q x_0\\
&+\int_0^t \mathbb E\left[x_{x_0}^T(s) Q x_{x_0}(s)\right] \left\|u(s)\right\|_2^2 ds,
\end{align*} 
where $y(t)=y(t, x_0, u)$. Due to Lemma \ref{gronwall}, we find \begin{align}\label{grobzglq}
\mathbb E\left[x^T_{x_0}(t) Q x_{x_0}(t)\right]\leq \alpha(t)+\int_0^t \alpha(s)\left\|u(s)\right\|_2^2 \exp\left(\int_s^t \left\|u(w)\right\|_2^2dw\right) ds,
\end{align} 
where we define $\alpha(t):=-\mathbb E\int_0^t \left\|y(s)\right\|_2^2ds+\mathbb E\int_0^t 2 x^T_{x_0}(s) Q B u(s)ds+x_0^T Q x_0$. We analyze this inequality further by looking at the terms depending on $x_0$:
\begin{align*}
&x_0^T Q x_0 \int_0^t \left\|u(s)\right\|_2^2 \exp\left(\int_s^t \left\|u(w)\right\|_2^2dw\right) ds= x_0^T Q x_0 \left[-\exp\left(\int_s^t \left\|u(w)\right\|_2^2dw\right)\right]_{s=0}^t\\
&=x_0^T Q x_0\left(\exp\left(\int_0^t \left\|u(s)\right\|_2^2ds\right)-1\right).
\end{align*}
Using this computation and $\mathbb E\left[x^T_{x_0}(t) Q x_{x_0}(t)\right]\geq 0$, we get from (\ref{grobzglq}) that \begin{align}\label{completeobserveq}
\mathbb E\int_0^t \left\|y(s)\right\|_2^2 ds \leq &x_0^T Q x_0 \exp\left(\int_0^t \left\|u(s)\right\|_2^2ds\right)\\
&+f_B(t)+\int_0^t f_B(s)\left\|u(s)\right\|_2^2 \exp\left(\int_s^t \left\|u(w)\right\|_2^2dw\right) ds,\nonumber
\end{align} 
where the term depending on the input matrix $B$ is $f_B(t):=\mathbb E\int_0^t 2 x^T_{x_0}(s) Q B u(s)ds$. 
In an observation problem an unknown initial condition $x_0$ is aimed to be reconstructed from the observations $y(t, x_0, u)$, $t\in[0, T]$. Since the control part $Bu$ does not depend on 
the unknown initial state, it can be assumed to be known and hence be neglected in the considerations by setting $B=0$. 
This assumption is also taken in \cite{bennerdamm, graymesko}, where the observation energy of deterministic bilinear systems is studied. 
Now, $B=0$ implies $f_B\equiv 0$. Applying this to (\ref{completeobserveq}) leads to the following bound on the observation energy on $[0, T]$: \begin{align}\label{bzeroobserveq}
\left.\mathbb E\int_0^T \left\|y(s, x_0, u)\right\|_2^2 ds \right\vert_{B=0} \leq x_0^T Q x_0 \exp\left(\int_0^T \left\|u(s)\right\|_2^2ds\right).
\end{align} 
If we the energy of the control is sufficiently small, e.g., $\left\|u\right\|_{L^2_T}\leq 1$, we can identify states from (\ref{bzeroobserveq}) producing only little observation energy.
We see that the energy that is caused by the observations of $x_0$ is small if the initial state is close to the kernel of $Q$. These initial states are contained in the eigenspaces of $Q$ corresponding to the small eigenvalues.
\section{$L^2$-error bound for BT}\label{errorboundsBT}

Let us assume that system (\ref{controlsystemoriginal}) has a zero initial condition ($x_0=0$) and is already balanced. Thus, (\ref{newgram2}) and (\ref{gengenlyapobs}) become
\begin{align}\label{balancedreach}
 A^T \Sigma^{-1}+\Sigma^{-1}A+\sum_{k=1}^m N_k^T \Sigma^{-1} N_k + \sum_{i, j=1}^v H_i^T \Sigma^{-1} H_j k_{i j} &\leq -\Sigma^{-1}BB^T \Sigma^{-1},\\ \label{balancedobserve}
 A^T \Sigma+\Sigma A+\sum_{k=1}^m N_k^T \Sigma N_k +\sum_{i, j=1}^v H_i^T \Sigma H_j k_{ij} &\leq -C^T C,
                                       \end{align}
i.e., $P=Q=\Sigma=\diag(\sigma_1, \ldots, \sigma_n)>0$. We partition the balanced coefficients of (\ref{controlsystemoriginal}) as follows:
\begin{align*}
 A=\smat{A}_{11}&{A}_{12}\\ 
{A}_{21}&{A}_{22}\srix,\;B=\smat B_1 \\ B_2\srix,\; N_k=\smat{N}_{k, 11}&{N}_{k, 12}\\ 
{N}_{k, 21}&{N}_{k, 22}\srix,\;H_i=\smat{H}_{i, 11}&{H}_{i, 12}\\ 
{H}_{i, 21}&{H}_{i, 22}\srix,\;C= \smat C_1 & C_2\srix,
            \end{align*}
where $A_{11}, N_{k, 11}, H_{i, 11}\in \mathbb R^{r\times r}$  ($k\in\left\{1, \ldots, m\right\}$ and $i\in\left\{1, \ldots, v\right\}$), $B_1\in \mathbb R^{r\times m}$ and $C_1\in \mathbb R^{p\times r}$ etc.
Furthermore, we partition the state variable and the Gramian \begin{align*}
               x=\smat x_1 \\ x_2\srix\text{ and }\Sigma=\smat \Sigma_1& \\ & \Sigma_2\srix,
                         \end{align*}
where $x_1$ takes values in $\mathbb R^r$ ($x_2$ accordingly), $\Sigma_1$ is the diagonal matrix of large HSVs and $\Sigma_2$ contains the small ones. The reduced system by BT is 
\begin{subequations}\label{romstochstatebt}
\begin{align}\label{romstateeq}
             dx_r&=[A_{11}x_r+B_1u+\sum_{k=1}^m N_{k, 11} x_r u_k]dt+\sum_{i=1}^v H_{i, 11} x_r dM_i,\\ 
    y_r(t)&=C_1x_r(t), \;\;\;t\geq 0,
            \end{align}
            \end{subequations}
where $x_r(0)=0$ and the time dependence in (\ref{romstateeq}) is omitted to shorten the notation. In order to find a bound for the approximation through BT, we define \begin{align*}
               x_-=\smat x_1-x_r \\ x_2\srix\text{ and } x_+=\smat x_1+x_r \\ x_2\srix,            
                           \end{align*}
and write down the corresponding equations for these variables. The system for $x_-$ is given by  
\begin{subequations}\label{xminus}
\begin{align}\label{statexminus}
             dx_-&=[Ax_-+\sum_{k=1}^m N_{k} x_- u_k]dt + \smat 0 \\ c_0\srix dt+\sum_{i=1}^v [H_{i} x_- + \smat 0 \\ c_i\srix] dM_i,\\ \label{xminusoutput}
    y_-(t)&=Cx_-(t)=Cx(t)-C_1 x_r(t)=y(t)-y_r(t), \;\;\;t\geq 0,
            \end{align}
            \end{subequations}
where $c_0(t):=A_{21}x_r(t)+B_2 u(t)+ \sum_{k=1}^m N_{k, 21} x_r(t) u_k(t)$ and $c_i(t):= H_{i, 21} x_r(t)$ for $i=1, \ldots, v$. We derive (\ref{xminus}) by comparing the partitioned system  (\ref{controlsystemoriginal}) 
with the reduced system (\ref{romstochstatebt}). The equation for $x_+$ looks similarly, the difference lies only in the signs for the compensation terms $c_0, \ldots, c_v$ and an additional control term:
\begin{align}\label{xplus}
             dx_+=[Ax_++2 B u+\sum_{k=1}^m N_{k} x_+ u_k]dt - \smat 0 \\ c_0\srix dt+\sum_{i=1}^v [H_{i} x_+ - \smat 0 \\ c_i\srix] dM_i.
            \end{align}
We will see that the proof of the error bound can be reduced to the task of finding suitable estimates for $\mathbb E[x_-^T(t) \Sigma x_-(t)]$ and $\mathbb E[x_+^T(t) \Sigma^{-1} x_+(t)]$.
The next theorem is the main result of this paper.
  \begin{thm}\label{mainthm}
Let $y$ be the output of the full model (\ref{controlsystemoriginal}) with $x(0)=0$ and $y_r$ be the output of the ROM (\ref{romstochstatebt}) with $x_{r}(0)=0$. Then, for all $T>0$, 
it holds that \begin{align*}
  \left(\mathbb E\left\|y-y_r\right\|_{L^2_{T}}^2\right)^{\frac{1}{2}}\leq 2 (\tilde \sigma_{1}+\tilde \sigma_{2}+\ldots + \tilde \sigma_\kappa)  \left\|u\right\|_{L^2_T}\exp\left(0.5 \left\|u\right\|_{L^2_T}^2\right),            
  \end{align*}
where $\tilde\sigma_{1}, \tilde\sigma_{2}, \ldots,\tilde\sigma_\kappa$ are the distinct diagonal entries of $\Sigma_2=\diag(\sigma_{r+1},\ldots,\sigma_n)=\diag(\tilde\sigma_{1} I, \tilde\sigma_{2} I, \ldots, \tilde\sigma_\kappa I)$. 
\begin{proof}
We compute an upper bound for $\mathbb E[x_-^T(t) \Sigma x_-(t)]$ making use of Lemma \ref{lemstochdiff}. Taking (\ref{statexminus}) into account then yields
\begin{align}\label{productruleapplied}
\mathbb E\left[x_-^T(t)\Sigma x_-(t)\right]=&2 \int_0^t\mathbb E\left[x_-^T\Sigma\left(Ax_-+\sum_{k=1}^m (N_{k} x_- u_k) + \smat 0 \\ c_0\srix \right)\right]ds\\ \nonumber
&+ \int_0^t\sum_{i, j=1}^v \mathbb E\left[\left(H_{i} x_- + \smat 0 \\ c_i\srix\right)^T\Sigma\left(H_{j} x_- + \smat 0 \\ c_j\srix\right)\right]k_{ij} ds,
\end{align}
where the time dependence of all functions is omitted for a shorter notation. Applying an estimate as in (\ref{maakebilinaerstoch}) provides
\begin{align*}
\sum_{k=1}^m 2 x_-^T(s)\Sigma N_{k} x_-(s) u_k(s) \leq  x_-^T(s) \Sigma x_-(s)  \left\|u(s)\right\|_{2}^2 +\sum_{k=1}^m  x_-^T(s) N_k^T\Sigma N_{k} x_-(s).
\end{align*}
Hence, (\ref{productruleapplied}) becomes 
\begin{align}\nonumber
\mathbb E\left[x_-^T(t)\Sigma x_-(t)\right]\leq& \mathbb E\int_0^t x_-^T\left(A^T \Sigma+\Sigma A+\sum_{k=1}^m  N_k^T\Sigma N_{k}+\sum_{i, j=1}^v H^T_{i}\Sigma H_{j}k_{ij}\right)x_- ds\\ \label{insetobserveequation}
& +\mathbb E\int_0^t 2 x_-^T\Sigma \smat 0 \\ c_0\srix + \sum_{i, j=1}^v \left(2 H_{i} x_- + \smat 0 \\ c_i\srix\right)^T\Sigma\smat 0 \\ c_j\srix k_{ij}  ds\\
&+\int_0^t \mathbb E\left[x_-^T \Sigma x_-\right]  \left\|u\right\|_{2}^2 ds.\nonumber
\end{align}
Using the partition for $x_-$ and $\Sigma$, we see that $x_-^T\Sigma \smat 0 \\ c_0\srix=x_2^T \Sigma_2 c_0$. With the partition of $H_i$, we additionally obtain 
\begin{align*}
&\left(2 H_{i} x_- + \smat 0 \\ c_i\srix\right)^T\Sigma\smat 0 \\ c_j\srix =\left(2 H_{i} x_- + \smat 0 \\ c_i\srix\right)^T\smat 0 \\ \Sigma_2 c_j\srix\\
&=\left(2 \smat H_{i, 21} & H_{i, 22}\srix (x - \smat x_r \\ 0\srix) +  c_i\right)^T \Sigma_2 c_j = \left(2 \smat H_{i, 21} & H_{i, 22}\srix x -  c_i\right)^T \Sigma_2 c_j.
\end{align*}
Inserting (\ref{balancedobserve}) and (\ref{xminusoutput}) into inequality (\ref{insetobserveequation}) and taking the above rearrangements into account leads to 
\begin{align*}
\mathbb E\left[x_-^T(t)\Sigma x_-(t)\right]\leq& - \mathbb E\left\|y-y_r\right\|^2_{L^2_{t}}+\int_0^t \mathbb E\left[x_-^T \Sigma x_-\right]  \left\|u\right\|_{2}^2 ds\\ 
& +\mathbb E\int_0^t 2 x_2^T \Sigma_2 c_0 + \sum_{i, j=1}^v \left(2 \smat H_{i, 21} & H_{i, 22}\srix x -  c_i\right)^T \Sigma_2 c_j k_{ij}  ds.
\end{align*}
We define $\alpha_-(t):=\mathbb E\int_0^t 2 x_2^T \Sigma_2 c_0 + \sum_{i, j=1}^v \left(2 \smat H_{i, 21} & H_{i, 22}\srix x -  c_i\right)^T \Sigma_2 c_j k_{ij}  ds$. Then, Lemma \ref{gronwall} implies
\begin{align*}
\mathbb E\left[x_-^T(t)\Sigma x_-(t)\right]\leq& \alpha_-(t)- \mathbb E\left\|y-y_r\right\|_{L^2_{t}}^2\\
& +\int_0^t (\alpha_-(s) - \mathbb E\left\|y-y_r\right\|_{L^2_{s}}^2) \left\|u(s)\right\|_{2}^2 \exp\left(\int_s^t \left\|u(w)\right\|_{2}^2 dw\right) ds.
\end{align*}
Thus, we find \begin{align*}
\mathbb E\left\|y-y_r\right\|_{L^2_{t}}^2\leq \alpha_-(t) +\int_0^t \alpha_-(s) \left\|u(s)\right\|_{2}^2 \exp\left(\int_s^t \left\|u(w)\right\|_{2}^2 dw\right) ds.
\end{align*}
We assume for the moment that $\Sigma_2=\sigma I$ and set $\alpha_+(t):=\mathbb E\int_0^t 2 x_2^T \Sigma_2^{-1} c_0 + \sum_{i, j=1}^v \left(2 \smat H_{i, 21} & H_{i, 22}\srix x -  c_i\right)^T \Sigma_2^{-1} c_j k_{ij} ds$. 
Hence, \begin{align}\label{firstbound}
\mathbb E\left\|y-y_r\right\|_{L^2_{t}}^2\leq \sigma^2\left[\alpha_+(t) +\int_0^t \alpha_+(s) \left\|u(s)\right\|_{2}^2 \exp\left(\int_s^t \left\|u(w)\right\|_{2}^2 dw\right) ds\right].
\end{align}
Let us turn our attention to the expression $\mathbb E[x_+^T(t) \Sigma^{-1} x_+(t)]$ for the further analysis of (\ref{firstbound}). Due to (\ref{xplus}) and Lemma \ref{lemstochdiff} it holds that
\begin{align}\label{productruleappliedplus}
\mathbb E\left[x_+^T(t)\Sigma^{-1} x_+(t)\right]=&2 \int_0^t\mathbb E\left[x_+^T\Sigma^{-1}\left(Ax_++2 Bu+\sum_{k=1}^m (N_{k} x_+ u_k) - \smat 0 \\ c_0\srix \right)\right]ds\\ \nonumber
&+ \int_0^t\sum_{i, j=1}^v \mathbb E\left[\left(H_{i} x_+ - \smat 0 \\ c_i\srix\right)^T\Sigma^{-1}\left(H_{j} x_+ - \smat 0 \\ c_j\srix\right)\right]k_{ij} ds.
\end{align}
As above, for the case of $x_-$, we use the estimate \begin{align*}
\sum_{k=1}^m 2 x_+^T(s)\Sigma^{-1} N_{k} x_+(s) u_k(s) \leq  x_+^T(s) \Sigma^{-1} x_+(s)  \left\|u(s)\right\|_{2}^2 +\sum_{k=1}^m  x_+^T(s) N_k^T\Sigma^{-1} N_{k} x_+(s),
\end{align*}
which leads to \begin{align}\nonumber
&\mathbb E\left[x_+^T(t)\Sigma^{-1} x_+(t)\right]\\ \nonumber&\leq \mathbb E\int_0^t x_+^T\left(A^T \Sigma^{-1}+\Sigma^{-1} A+\sum_{k=1}^m  N_k^T\Sigma^{-1} N_{k}+\sum_{i, j=1}^v H^T_{i}\Sigma^{-1} H_{j}k_{ij}\right)x_+ ds\\ \label{insetobserveequation2}
&\quad -\mathbb E\int_0^t 2 x_+^T\Sigma^{-1} \smat 0 \\ c_0\srix + \sum_{i, j=1}^v \left(2 H_{i} x_+ - \smat 0 \\ c_i\srix\right)^T\Sigma^{-1}\smat 0 \\ c_j\srix k_{ij}  ds\\
&\quad+\int_0^t \mathbb E\left[x_+^T \Sigma^{-1} x_+\right]  \left\|u\right\|_{2}^2 ds+ 4 \mathbb E \int_0^t x_+^T\Sigma^{-1} Bu ds.\nonumber
\end{align}
From inequality (\ref{balancedreach}) and the Schur complement condition on definiteness it follows that\begin{align}\label{schurposdef}
 \mat{cc}\hspace{-0.15cm} A^T \Sigma^{-1}\hspace{-0.025cm}+\Sigma^{-1}\hspace{-0.025cm}A+\hspace{-0.025cm}\sum_{k=1}^m N_k^T \Sigma^{-1} N_k + \hspace{-0.025cm}\sum_{i, j=1}^v H_i^T \Sigma^{-1} H_j k_{i j}  & \Sigma^{-1}B\\
 B^T \Sigma^{-1}& -I\rix\leq 0.
                                       \end{align}
We multiply (\ref{schurposdef}) with $\smat x_+ \\ 2u\srix^T$ from the left and with $\smat x_+\\ 2u\srix$ from the right. This leads to 
\begin{align*}
&4 \left\|u\right\|_{2}^2\geq \\
& x_+^T\left(A^T \Sigma^{-1}+\Sigma^{-1} A+\sum_{k=1}^m  N_k^T\Sigma^{-1} N_{k}+\sum_{i, j=1}^v H^T_{i}\Sigma^{-1} H_{j}k_{ij}\right)x_++4x_+^T\Sigma^{-1} Bu.
\end{align*}
Applying this result to inequality (\ref{insetobserveequation2}) gives 
\begin{align}\label{insetobserveequationbla}
\mathbb E\left[x_+^T(t)\Sigma^{-1} x_+(t)\right] \leq & 4 \left\|u\right\|_{L^2_t}^2+\int_0^t \mathbb E\left[x_+^T \Sigma^{-1} x_+\right]  \left\|u\right\|_{2}^2 ds\\ \nonumber
& -\mathbb E\int_0^t 2 x_+^T\Sigma^{-1} \smat 0 \\ c_0\srix + \sum_{i, j=1}^v \left(2 H_{i} x_+ - \smat 0 \\ c_i\srix\right)^T\Sigma^{-1}\smat 0 \\ c_j\srix k_{ij}  ds.
\end{align}
We further analyze the terms in (\ref{insetobserveequationbla}). We find that $x_+^T\Sigma^{-1} \smat 0 \\ c_0\srix=x_2^T\Sigma_2^{-1}c_0$ using the partitions of $x_+$ and $\Sigma$. 
 With the partition of $H_i$, we moreover have 
\begin{align*}
&\left(2 H_{i} x_+ - \smat 0 \\ c_i\srix\right)^T\Sigma^{-1}\smat 0 \\ c_j\srix =\left(2 H_{i} x_+ - \smat 0 \\ c_i\srix\right)^T\smat 0 \\ \Sigma_2^{-1} c_j\srix\\
&=\left(2 \smat H_{i, 21} & H_{i, 22}\srix (x + \smat x_r \\ 0\srix) -  c_i\right)^T \Sigma_2^{-1} c_j = \left(2 \smat H_{i, 21} & H_{i, 22}\srix x +  c_i\right)^T \Sigma_2 c_j.
\end{align*}
We plug this into (\ref{insetobserveequationbla}), such that 
\begin{align}\label{insetobserveequationblabla}
\mathbb E\left[x_+^T(t)\Sigma^{-1} x_+(t)\right] \leq & 4\left\|u\right\|_{L^2_t}^2+\int_0^t \mathbb E\left[x_+^T \Sigma^{-1} x_+\right]  \left\|u\right\|_{2}^2 ds\\ \nonumber
& -\mathbb E\int_0^t 2 x_2^T\Sigma_2^{-1} c_0 + \sum_{i, j=1}^v \left(2 \smat H_{i, 21} & H_{i, 22}\srix x + c_i\right)^T\Sigma_2^{-1} c_j k_{ij} ds.
\end{align}
Adding $2\mathbb E\int_0^t \sum_{i, j=1}^v c_i^T\Sigma_2^{-1} c_j k_{ij} ds$ to the right side of (\ref{insetobserveequationblabla}), which is a nonnegative term due to Lemma \ref{proppossemidef}, we have
\begin{align*}
\mathbb E\left[x_+^T(t)\Sigma^{-1} x_+(t)\right] \leq 4 \left\|u\right\|_{L^2_t}^2-\alpha_+(t)+\int_0^t \mathbb E\left[x_+^T(s) \Sigma^{-1} x_+(s)\right]  \left\|u(s)\right\|_{2}^2 ds.
\end{align*}
Lemma \ref{gronwall} yields \begin{align}\label{keineahnung}
\mathbb E\left[x_+^T(t)\Sigma^{-1} x_+(t)\right] \leq & 4\left\|u\right\|_{L^2_t}^2-\alpha_+(t)\\ \nonumber
&+\int_0^t (4 \left\|u\right\|_{L^2_s}^2-\alpha_+(s))  \left\|u(s)\right\|_{2}^2 \exp\left(\int_s^t \left\|u(w)\right\|_{2}^2 dw\right) ds.
\end{align}
Moreover, we have \begin{align}\nonumber
&\int_0^t \left\|u\right\|_{L^2_s}^2 \left\|u(s)\right\|_2^2 \exp\left(\int_s^t \left\|u(w)\right\|_2^2dw\right) ds\leq \left\|u\right\|_{L^2_t}^2 \left[-\exp\left(\int_s^t \left\|u(w)\right\|_2^2dw\right)\right]_{s=0}^t\\
&=\left\|u\right\|_{L^2_t}^2\left(\exp\left(\int_0^t \left\|u(s)\right\|_2^2ds\right)-1\right). \label{keineahnung2}
\end{align}
Combining (\ref{keineahnung}) with (\ref{keineahnung2}), we get
\begin{align*}
\alpha_+(t)+\int_0^t \alpha_+(s)  \left\|u(s)\right\|_{2}^2 \exp\left(\int_s^t \left\|u(w)\right\|_{2}^2 dw\right) ds\leq 4\left\|u\right\|_{L^2_t}^2 \exp\left(\int_0^t \left\|u(s)\right\|_2^2ds\right).
\end{align*}
Comparing this result with (\ref{firstbound}) implies 
 \begin{align}\label{secondbound}
\left(\mathbb E\left\|y-y_r\right\|_{L^2_{t}}^2\right)^{\frac{1}{2}}\leq 2\sigma \left\|u\right\|_{L^2_t} \exp\left(0.5 \left\|u\right\|_{L^2_t}^2\right).
\end{align}
For the proof of a general $\Sigma_2$, we remove the HSVs step by step. We use the triangle inequality to bound the error between the outputs $y$ and $y_r$:\begin{align*}
  &\left(\mathbb E \left\|y-y_r\right\|^2_{L^2_T}\right)^{\frac{1}{2}}\\&\leq   
\left(\mathbb E\left\|y-y_{r_\kappa}\right\|^2_{L^2_T}\right)^{\frac{1}{2}}+\left(\mathbb E\left\|y_{r_\kappa}-y_{r_{\kappa-1}}\right\|^2_{L^2_T}\right)^{\frac{1}{2}}+\ldots+\left(\mathbb E\left\|y_{r_2}-y_{r}\right\|^2_{L^2_T}\right)^{\frac{1}{2}},        
  \end{align*}
where the dimensions $r_i$ of the corresponding states are defined by $r_{i+1}=r_{i}+m(\tilde\sigma_{i})$ for $i=1, 2 \ldots, \kappa-1$. The number 
$m(\tilde\sigma_{i})$ denotes the multiplicity of $\tilde\sigma_{i}$ and $r_1=r$. In the first step only the smallest 
HSV $\tilde\sigma_\kappa$ is removed from the system. By inequality (\ref{secondbound}), we have \begin{align*}             
  \left(\mathbb E\left\|y-y_{r_\kappa}\right\|^2_{L^2_T}\right)^{\frac{1}{2}}\leq 2 \tilde\sigma_\kappa \left\|u\right\|_{L^2_T} \exp\left(0.5 \left\|u\right\|_{L^2_T}^2\right).
\end{align*}
The same kind of bound can be established when comparing the reduced order outputs $y_{r_\kappa}$ and $y_{r_{\kappa-1}}$. Again, only 
one HSV, namely $\tilde\sigma_{r_{\kappa-1}}$, is removed. Moreover, the matrix inequalities in the ROM have the same structure as (\ref{balancedreach}) and (\ref{balancedobserve}). 
To be more precise, evaluating the left upper blocks of (\ref{balancedreach}) and (\ref{balancedobserve}), we obtain\begin{align*}
 A_{11}^T \Sigma_1^{-1}+\Sigma_1^{-1}A_{11}+\sum_{k=1}^m N_{k, 11}^T \Sigma_1^{-1} N_{k, 11}+\sum_{i, j=1}^v H_{i, 11}^T \Sigma_1^{-1} H_{j, 11} k_{ij} &\leq -\Sigma_1^{-1}B_1B_1^T
\Sigma_1^{-1},\\
 A_{11}^T \Sigma_1+\Sigma_1 A_{11}+\sum_{k=1}^m N_{k, 11}^T \Sigma_1 N_{k, 11}+\sum_{i, j=1}^v H_{i, 11}^T \Sigma_1 H_{j, 11} k_{ij} &\leq -C_1^TC_1
                                       \end{align*}
applying Lemma \ref{proppossemidef}. Thus, by repeatedly applying the above arguments, we have\begin{align*}
  \left(\mathbb E \left\|y_{r_j}-y_{r_{j-1}}\right\|^2_{L^2_T}\right)^{\frac{1}{2}}\leq 2\tilde\sigma_{r_{j-1}} \left\|u\right\|_{L^2_T} \exp\left(0.5 \left\|u\right\|_{L^2_T}^2\right)
  \end{align*}
for $j=2, \ldots, \kappa $. This concludes the proof.
\end{proof}
\end{thm}\\
The result in Theorem \ref{mainthm} is the first one of that type for deterministic/stochastic systems. In contrast to \cite{redmanntypeiibilinear}, we assume no bound on the control $u$ which possibly can be small. 
We pay a price for dealing with general $L_T^2$ controls, since we obtain an exponential term in Theorem \ref{mainthm} which is due to the bilineararity. However, this result is still very meaningful, because it tells us that 
the ROM (\ref{romstochstatebt}) yields a very good approximation if the truncated HSVs (diagonal entries of $\Sigma_2$) are small and, e.g., a normalized control is used. At the same time, the exponential term 
in the error bound can be an indicator that BT performs terribly bad if the control energy is large.
\section{Conclusions}
In this paper, we investigated a large-scale stochastic bilinear system. In order to reduce the state space dimension, a model order reduction technique called balanced truncation was extended to this setting. 
To do so, we proposed a reachability and an observability Gramian. We proved energy estimates with the help of these Gramian that allow us to find the unimportant states within the system. The reduced system 
was then obtained by removing these states from the stochastic bilinear system. Finally, we provided a new error bound that can be used to point out the cases in which the reduced order model by balanced truncation 
delivers a good approximation to the original model.
\appendix

\section{Supporting Lemmas}

In this appendix, we state three important results and the corresponding references that we frequently use throughout this paper.
\begin{lem}\label{lemstochdiff}
Let $a, b_1, \ldots, b_v$ be $\mathbb R^d$-valued processes, where $a$ is adapted and almost surely Lebesgue integrable and the functions $b_i$ 
are integrable with respect to the mean zero square integrable L\'evy process $M=(M_1, \ldots, M_v)^T$ with covariance matrix $K=\left(k_{ij}\right)_{i, j=1, \ldots, v}$. If the process $x$ is given by \begin{align*}
 dx(t)=a(t) dt+ \sum_{i=1}^v b_i(t)dM_i,                                                                  
                                                                   \end{align*}
then, we have \begin{align*}
 \frac{d}{dt}\mathbb E\left[x^T(t) x(t)\right]=2 \mathbb E\left[x^T(t) a(t)\right] + \sum_{i, j=1}^v \mathbb E\left[b_i^T(t) b_j(t)\right]k_{ij}.  
                                                                                                                                   \end{align*}
\begin{proof}
We refer to \cite[Lemma 5.2]{redmannspa2} for a proof of this lemma. 
\end{proof}
\end{lem}
\begin{lem}\label{proppossemidef}
Let $A_1, \ldots, A_v$ be $d_1\times d_2$ matrices and $K=(k_{ij})_{i, j=1, \ldots, v}$ be a positive semidefinite matrix, then
\begin{align*}\tilde K:=\sum_{i,j=1}^v A_i^T A_j k_{ij} \end{align*}
is also positive semidefinite.
 \begin{proof}
  The proof can be found in \cite[Proposition 5.3]{redmannspa2}.
 \end{proof}
\end{lem}

\begin{lem}[Gronwall lemma]\label{gronwall}
Let $T>0$, $z, \alpha: [0, T]\rightarrow \mathbb R$ be measurable bounded functions and $\beta: [0, T]\rightarrow \mathbb R$ be a nonnegative integrable function. 
If \begin{align*}
    z(t)\leq \alpha(t)+\int_0^t \beta(s) z(s) ds,
   \end{align*}
then it holds that \begin{align}\label{gronwallineq}
    z(t)\leq \alpha(t)+\int_0^t \alpha(s)\beta(s) \exp\left(\int_s^t \beta(w)dw\right) ds
   \end{align}
for all $t\in[0, T]$.  Moreover, a non-decreasing function $\alpha$ implies that (\ref{gronwallineq}) becomes \begin{align*}
    z(t)\leq \alpha(t)\exp\left(\int_0^t \beta(s)ds\right)
   \end{align*}
   for all $t\in[0, T]$.
 \begin{proof}
  The result is shown as in \cite[Proposition 2.1]{gronwalllemma}.
 \end{proof}
\end{lem}

\bibliographystyle{abbrv}

\end{document}